\documentclass[opre,nonblindrev]{informs4}
\OneAndAHalfSpacedXI

\usepackage{natbib}
 \bibpunct[, ]{(}{)}{,}{a}{}{,}%
 \def\bibfont{\small}%
 \def\bibsep{\smallskipamount}%
\usepackage{longtable}

\usepackage{xr}
\makeatletter
\newcommand*{\addFileDependency}[1]{
  \typeout{(#1)}
  \@addtofilelist{#1}
  \IfFileExists{#1}{}{\typeout{No file #1.}}
}
\makeatother

\usepackage{setspace}
\usepackage{url}
\usepackage{booktabs}
\usepackage{amsfonts}
\usepackage{nicefrac}
\usepackage{microtype}
\usepackage{float}
\usepackage{enumerate}
\usepackage{enumitem}
\usepackage{listings}
\usepackage{mathtools}
\usepackage{amsmath}
\usepackage{amssymb}
\usepackage{centernot}
\usepackage{float}
\usepackage{braket}
\usepackage{wrapfig}

\TheoremsNumberedThrough
\ECRepeatTheorems
\EquationsNumberedThrough

\begin{document}
%%%%%%%%%%%%%%%%
\MONTH{September}
\YEAR{2020}
\RUNTITLE{Diet Dataset for Constrained Optimization Models}
\ARTICLEAUTHORS{%
\AUTHOR{Farzin Ahmadi, Fardin Ganjkhanloo, Kimia Ghobadi}
\AFF{Department of Civil and Systems Engineering, The Center for Systems Science and Engineering, The Malone Center for Engineering in Healthcare, Johns Hopkins University, Baltimore, MD 21218, \EMAIL{fahmadi1@jhu.edu}, \EMAIL{fganjkh1@jhu.edu}, \EMAIL{kimia@jhu.edu}} %, \URL{}} 
% Enter all authors
}

\TITLE{An Open-Source Dataset on Dietary Behaviors \\  and DASH Eating Plan Optimization Constraints}

\RUNAUTHOR{Ahmadi, Ganjkhanloo, Ghobadi}
\ABSTRACT{%
Linear constrained optimization techniques have been applied to many real-world settings. In recent years, inferring the unknown parameters and functions inside an optimization model has also gained traction. This inference is often based on 
%there has been an increase in attention towards inference of unknown parameters and functions inside an optimization model based on 
existing observations and/or known parameters. Consequently, %many of 
such models require reliable, easily accessed, and easily interpreted examples to be evaluated. 
To facilitate research in such directions, we provide a modified dataset based on dietary behaviors of different groups of people, their demographics, and pre-existing conditions, among other factors. This data is gathered from the National Health and Nutrition Examination Survey (NHANES) and complemented with the nutritional data from the United States Department of Agriculture (USDA). We additionally provide tailored datasets for hypertension and pre-diabetic patients as groups of interest who may benefit from targetted diets such as the Dietary Approaches to Stop Hypertension (DASH) eating plan. The data is compiled and curated in such a way that it is suitable as input to linear optimization models. 
We hope that this data and its supplementary, open-accessed materials can accelerate and simplify interpretations and research on linear optimization and constrained inference models. The complete dataset can be found in the following repository: https://github.com/CSSEHealthcare/InverseLearning
}

\maketitle
\section{Introduction}

Data-driven approaches in constrained inference models and constrained optimization problems are increasingly gaining traction. Considering their immense impact on applied settings, it deems necessary to provide researchers with open-access and accurate real-world examples and datasets to be utilized in research and methodology evaluation. To this end, we provide a dataset of dietary behavior within the diet recommendation setting. This dataset is gathered and curated from multiple sources and provides researchers with an interpretable and accurate (partially self-reporting) data that can be readily used in models and approaches. 

The classical diet recommendation optimization problem was first posed by \cite{stigler1945cost} as finding an optimal set of intakes of different foods with minimum cost while satisfying a given set of minimum nutrient conditions for an average person \citep{stigler1945cost,garille2001stigler}. We consider this renowned optimization problem in practical settings, where diet plays a vital role in controlling the progression of chronic illnesses such as hypertension and type II diabetes. We provide specific nutritional and food group constraints for a typical diet that clinicians recommend for hypertension patients, the Dietary Approaches to Stop Hypertension (DASH) eating plan \citep{sacks2001effects,liese2009adherence}. In the rest of this technical report, we provide details about our data collection and insights into its use in optimization and inference models. 

\section{Dataset}

This diet behavior dataset includes all data points necessary to perform inference analysis in a data-driven optimization setting. Specifically, the dataset is apt for studying dietary behaviors in %a 
diverse demographics or providing personalized recommendations in constrained environments. 
Our dataset includes two main components of (1) individuals' daily consumption details and (2) the dietary requirements. The first component (details of individuals' diets) contains information gathered from the National Health and Nutrition Examination Survey \citep{CDC_2020} and is augmented with the United States Department of Agriculture (USDA) data for food nutrients.   
The second component of our dataset focuses on the DASH diet as a potential target diet to construct nutrient constraints. Throughout our dataset, we use servings as the unit of measurement in our data to provide a more general ground for comparison. 

{\bf Individuals' daily consumption.}  
We extracted a set of daily dietary habits of individuals from the NHANES Dietary data. This data is self-reported information on the daily intake of 9,544 patients over two days. This input data from NHANES is divided into two different sets. The first set contains two days of detailed food intake data includes the type, the amount, and the nutrients of each food. The second dataset is the aggregate nutrient intakes of each respondent for the two days of measurement. We further group the patients based on their demographics and whether or not they are concerned by or diagnosed with hypertension.

The intake data include more than 5,000 different food types. Given the large number of food types, we bundled them into 38 broad food groups for ease of interpretations and to make the inferred diets more tractable. This categorization is done based on the food codes from USDA.
Table \ref{Table:food_groups} shows the grouping developed for the dataset and the average serving size in grams. We note that although the introduction of such food categorization affects the accuracy of the results of the models, the benefits gained in the size reductions greatly outweigh the approximations. 

\begin{table}[] 
\small
        \caption{Food groups and their respective serving sizes in grams}
        \label{Table:food_groups}
\begin{tabular}{>{}p{0.24\textwidth}|>{}p{0.55\textwidth}|>{\centering\arraybackslash}p{0.15\textwidth}}
Group Name              & Description                                                                                            & Serving Size (g) \\ \hline \hline
Milk                    & milk, soy milk, almond milk, chocolate milk, yogurt, baby food, infant   formula                         & 244              \\
Cream                   & Cream, sour cream                                                                                      & 32               \\
Ice Cream               & all types of ice cream                                                                                 & 130              \\
Cheese                  & all types of cheese                                                                                    & 32               \\
Beef                    & ground beef, steaks (cooked, boiled, grilled or raw)                                                   & 65               \\
Pork                    & chops of pork, cured pork, bacon (cooked, boiled, grilled or raw)                                      & 84               \\
Red Meat (Other)        & lamb, goat, veal, venison  (cooked, boiled, grilled or raw)                                            & 85               \\
Chicken, Turkey          & all types of chicken, turkey, duck (cooked, boiled, grilled or raw)                                   & 110              \\
Sausages                & beef or red meat by-products, bologna, sausages, salami, ham    (cooked, boiled, grilled or raw)      & 100              \\
Fish                    & all types of fish,                                                                                     & 85               \\
Stew                    & stew meals containing meat (or substitutes), rice, vegetables                                          & 140              \\
Frozen Meals            & frozen meal (containing meat and vegetables)                                                           & 312              \\
Egg Meals               & egg meals, egg omelets and substitutes                                                                 & 50               \\
Beans                   & all types of beans (cooked, boiled, baked, raw)                                                                                     & 130           \\
Nuts                    & all types of nuts                                                                                      & 28.35            \\
Seeds                   & all types of seeds                                                                                     & 30               \\
Bread                   & all types of bread                                                                                     & 25               \\
Cakes, Biscuits, Pancakes & cakes, cookies, pies, pancakes, waffles                                                                & 56               \\
Noodle, Rice             & macaroni, noodle, pasta, rice                                                                          & 176              \\
Cereal                  & all types of cereals                                                                                   & 55               \\
Fast Foods              & burrito, taco, enchilada, pizza, lasagna                                                               & 198              \\
Meat Substitutes        & meat substitute that are  cereal-   or vegetable protein-based                                         & 100              \\
Citrus Fruits           & grapefruits, lemons, oranges                                                                           & 236              \\
Dried Fruits            & all types of dried fruit                                                                               & 28.35            \\
Tropical Fruits         & apples, apricots, avocados, bananas, cantaloupes, cherries, figs, grapes,   mangoes, pears, pineapples & 182              \\
Fruit Juice             & All types of fruit juice                                                                               & 249              \\
Potato products                & potatoes (fried, cooked)                                                                                           & 117              \\
Greens                  & beet greens, collards, cress, romaine, greens, spinach                                                 & 38               \\
Squash/Roots            & carrots, pumpkins, squash, sweet potatoes                                                               & 72               \\
Tomato products                 & tomato, salsa containing tomatoes, tomato byproducts                                                   & 123              \\
Vegetables              & raw vegetables                                                                                         & 120              \\
Puerto Rican Food      & Puerto Rican style food                                                                                & 250              \\
Smoothies               & fruit and vegetable smoothies                                                                          & 233              \\
Butter, Oils             & butters, oils                                                                                          & 14.2             \\
Salad Dressing          & all types of salad dressing                                                                            & 14               \\
Desserts                & sugars, desserts, toppings                                                                             & 200              \\
Caffeinated Drinks      & coffees, soda drinks, iced teas                                                                        & 240              \\
Shakes                  & shakes, drink mixes                                                                                    & 166   \\    
\hline
\end{tabular}
\end{table}

{\bf Dietary Requirements.} 
The DASH diet (Dietary Approaches to Stop Hypertension) is developed to lower blood pressure without medication \citep{mayo_clinic_2019} by limiting the amount of sodium intake of patients. Table \ref{Table:dash_diet_recommendations_servings} illustrates the recommendations of the DASH diet in terms of the number of servings of each food group for different diets with distinct calorie targets. Since the DASH diet recommendations are in servings, Table \ref{Table:dash_diet_recommendations_servings} provides additional details about a typical sample of each food group along with the corresponding amount in one serving size.  We utilize the food samples from Table \ref{Table:dash_diet_recommendations_servings}, the nutritional data from USDA, and the recommended amounts from the DASH eating plan to calculate the required bounds on nutrients. These bounds can serve as the right-hand side vector for constraints in linear optimization settings. Table \ref{Table:SubsetNutrients} depicts these amounts for the 1,600 and 2,000 calorie target diets. Additional nutritional bounds for other diets can be found in the repository provided for this work \citep{CSSEDietData}.

\begin{table}[]
\small
        \caption{Food categories and their recommended number of servings for different targets based on the DASH diet \citep{dash_diet_2020}}
        \label{Table:dash_diet_recommendations_servings}
\begin{tabular}{>{}p{0.3\textwidth}|>{}p{0.1\textwidth}|>{}p{0.08\textwidth}|>{}p{0.08\textwidth}|>{}p{0.08\textwidth}|>{}p{0.08\textwidth}|>{}p{0.075\textwidth}|>{\arraybackslash}p{0.075\textwidth}} 
&\multicolumn{7}{c}{Diet Target}\\\cline{2-8}
%&&&&Diet Calorie Target &&&\\
Food Category                          & 1,200    \ \  Calories               & 1,400 Calories              & 1,600 Calories               & 1,800 Calories               & 2,000 Calories               & 2,600 Calories         & 3,100 Calories\\
\hline \hline
Grains                             & 4–5                & 5–6                & 6                  & 6                  & 6–8                & 10–11        & 12–13        \\
Vegetables                          & 3–4                & 3–4                & 3–4                & 4–5                & 4–5                & 5–6          & 6            \\
Fruits                              & 3–4                & 4                  & 4                  & 4–5                & 4–5                & 5–6          & 6            \\
Fat-free or low-fat dairy products & 2–3                & 2–3                & 2–3                & 2–3                & 2–3                & 3            & 3–4          \\
Lean meats, poultry, and fish       & $\leq$ 3          & $\leq$3–4        & $\leq$3–4         & $\leq$ 6         & $\leq$6          & $\leq$6     & 6–9          \\
Nuts, seeds, and legumes            & 3/week         & 3/week         & 3–4/week       & 4/week         & 4–5/week       & 1            & 1            \\
Fats and oils                      & 1                  & 1                  & 2                  & 2–3                & 2–3                & 3            & 4            \\
Sweets and added sugars             & $\leq$ 3/week & $\leq$3/week & $\leq$3/week & $\leq$5/week & $\leq$5/week & $\leq$2           & $\leq$2           \\
Maximum sodium limit(mg/day)               & 2,300        & 2,300        & 2,300        & 2,300       & 2,300        & 2,300  & 2,300 
\end{tabular}
\end{table}

\begin{table}[]
\small
        \caption{Food categories and their respective serving sizes in grams}
        \label{Table:dash_diet_servings}
\begin{tabular}{ll}
Food Category                            & Serving Size (Example)                                                        \\
\hline \hline
Grains                              & 1 slice of whole-grain bread                                                 \\
Vegetables                            & 1 cup (about 30 grams) of raw, leafy green vegetables like   spinach or kale \\
Fruits                                & 1 medium apple                                                               \\
Fat-free   or low-fat dairy products & 1 cup (240 ml) of low-fat milk                                               \\
Lean   meats, poultry, and fish       & 1 ounce (28 grams) of cooked meat, chicken or fish                           \\
Nuts,   seeds, and legumes            & 1/3 cup (50 grams) of nuts                                                   \\
Fats   and oils                     & 1 teaspoon (5 ml) of vegetable oil                                           \\
Sweets   and added sugars             & 1 cup (240 ml) of lemonade                                                  
\end{tabular}
\end{table}

\begin{table}[]
    \caption{Nutrient bounds based on select diets from the recommendations of the DASH eating plans.}
    \label{Table:SubsetNutrients}
\begin{tabular}{>{}p{0.3\textwidth}|>{}p{0.15\textwidth}>{}p{0.15\textwidth}|>{}p{0.15\textwidth}>{\arraybackslash}p{0.15\textwidth}}
                          & \multicolumn{2}{c}{1,600 Calories} & \multicolumn{2}{c}{2,000 Calories} \\
Nutrients                 & Lower Bound   & Upper Bound   & Lower Bound   & Upper Bound   \\
\hline \hline
Food energy (kcal)        & 1,575.1        & 2,013.2        & 1,663.0        & 2,611.0        \\
Carbohydrate (g)          & 223.7         & 254.2         & 229.1         & 320.8         \\
Protein (g)               & 51.8          & 89.2          & 55.0          & 113.9         \\
Total Fat(g)              & 59.7          & 78.2          & 66.0          & 106.0         \\
Total Sugars (g)          & 117.0         & 144.6         & 119.5         & 177.1         \\
Fiber (g)                 & 36.7          & 39.3          & 39.3          & 50.1          \\
Saturated fatty acids (g) & 11.4          & 16.6          & 12.4          & 21.4          \\
Cholesterol (mg)          & 24.4          & 120.6         & 24.4          & 162.6         \\
Iron (mg)                 & 9.7           & 12.5          & 10.9          & 16.3          \\
Sodium (mg)               & 1,376.2        & 1,693.0        & 1,473.2        & 2,137.3        \\
Caffeine (mg)             & 0.0           & 80.0           & 0.0           & 80.0          
\end{tabular}
\end{table}

\section{Discussion}
With the increasing applications and importance of optimization and inference models, the ease of access to reliable, accurate, and interpretable datasets has become paramount. In this report, we provided the details of a large-scale open-access dataset on dietary behaviors and attributes of individuals. We complemented this dataset with data on nutritional requirements that are typical for hypertension and pre-diabetic patients. We hope that the presence of such a dataset can help the researchers in different data-driven approaches to evaluate proposed methods and get meaningful insights. 

\bibliographystyle{informs2014}
\bibliography{Tech_Report}

\end{document}